\documentclass[11pt]{article}
\usepackage{epsf,amsmath,amsfonts,amssymb,amsthm}
\usepackage{enumerate}

\theoremstyle{plain}						
\newtheorem{proposition}{Proposition}[section]
\newtheorem{lemma}[proposition]{Lemma}
\newtheorem{theorem}[proposition]{Theorem}
\newtheorem{corollary}[proposition]{Corollary}

\theoremstyle{definition}

\theoremstyle{remark}

\font\tenscr	= rsfs10 
\font\sevenscr	= rsfs7  
\font\fivescr	= rsfs5  

\font\tenbbm	= bbm10
\font\sevenbbm	= bbm7
\font\fivebbm	= bbm5

\newfam\bbmfam
\textfont\bbmfam=\tenbbm
\scriptfont\bbmfam=\sevenbbm
\scriptscriptfont\bbmfam=\fivebbm

\font\tenscr=rsfs10 
\font\sevenscr=rsfs7 
\font\fivescr=rsfs5 
\skewchar\tenscr=127 \skewchar\sevenscr=127 \skewchar\fivescr=127
\newfam\scrfam
\textfont\scrfam=\tenscr
\scriptfont\scrfam=\sevenscr
\scriptscriptfont\scrfam=\fivescr
\def\scr{\fam\scrfam}

\input scrload

\def\ccc{\mathfrak{c}} 		
\def\I {{\scr I}}		
\def\F {{\scr F}}		
\def\G {{\scr G}}		
\def\U {{\scr U}}		
\def\W {{\scr W}}		
\def\D {{\mathcal{D}}}	

\begin{document}

\begin{center}
{\LARGE The relation of rapid ultrafilters and Q-points\\ to van der Waerden ideal}\\[5mm]
\textsc{Jana Fla\v{s}kov\'a}
\footnote{The current research was partially supported by the Institut Mittag-Leffler (Djursholm, Sweden) and by the European Science Foundation (in the realm of the activity entitled 'New Frontiers of Infinity: Mathematical, Philosophical and Computational Prospects').}\\[3mm]
Department of Mathematics, University of West Bohemia
\end{center}


\bigskip

\begin{abstract}
We point out one of the differences between rapid ultrafilters and $Q$-points: Rapid ultrafilters
may have empty intersection with the van der Waerden ideal, whereas every $Q$-point has a non-empty
intersection with the van der Waerden ideal. Assuming Martin's axiom for countable posets we also construct
a $\W$-ultrafilter which is not a $Q$-point.
\end{abstract}

\section{Preliminaries}

Let us first recall the definitions of the two classes of ultrafilters in question:
An ultrafilter $\U$ is a {\it $Q$-point\/} if for every partition
$\{Q_n: n \in \omega\}$ of $\omega$ into finite sets there exists $A \in \U$
such that $|A \cap Q_n| \leq 1$ for every $n \in \omega$.\\
An ultrafilter $\U$ is a {\it rapid ultrafilter\/} if the
enumeration functions of sets in $\U$ form a dominating family in
$({}^{\omega}\omega, \leq^{\ast})$, where the enumeration function of a set $A$ is
the unique strictly increasing function $e_A$ from $\omega$ onto $A$.\\
Every $Q$-point is a rapid ultrafilter, but the converse is not true (see \cite{Mi}).
\smallskip

Remember that a set $A \subseteq \omega$ is called an {\it AP-set\/} if it contains arbitrary
long arithmetic progressions. It follows from the van der Waerden theorem that
sets which are not AP-sets form a proper ideal on $\omega$. We will refer to
this ideal as {\it van der Waerden ideal\/} and denote it by $\W$. It is known that
the van der Waerden ideal is an $F_{\sigma}$-ideal.\\
An ideal $\I$ on $\omega$ is {\it tall\/} if for every infinite $A \subseteq \omega$ there exists
infinite $B \subseteq A$ such that $B \in \I$. The van der Waerden ideal is clearly a tall ideal.\\
An ideal $\I$ on $\omega$ is a {\it P-ideal\/} if for every $A_n \in \I$, $n \in \omega$ there
exists $A \in \I$ such that $A_n \setminus A$ is finite for all $n \in \omega$. It is easy to see
that the van der Waerden ideal is not a $P$-ideal (consider e.g.~sets $A_n = \{2^i + n: i \in \omega\}$).
\smallskip

Baumgartner introduced $\I$-ultrafilters in \cite{B}. We will repeat the definition and
introduce a weaker notion of weak $\I$-ultrafilters. So, assume $\I$ is a tall ideal on $\omega$:\\
An ultrafilter $\U$ on $\omega$ is called an {\it $\I$-ultrafilter\/} if
for every function $f:\omega \rightarrow \omega$ there exists $U \in \U$
such that $f[U] \in \I$.\\
An ultrafilter $\U$ on $\omega$ is called {\it a weak $\I$-ultrafilter} if
for every finite-to-one function $f:\omega \rightarrow \omega$ there exists
$U \in \U$ such that $f[U] \in \I$.\\
It is obvious that every $\I$-ultrafilter is a weak $\I$-ultrafilter and every weak $\I$-ultrafilter
has nonempty intersection with the ideal $\I$.

\section{$\mathbf{Q}$-points and the van der Waerden ideal}

In this section we examine the connection between $Q$-points and the van der Waerden ideal.
First of all we prove that each $Q$-point not only has a nonempty intersection with the van der Waerden
ideal, but it is in fact a weak $\W$-ultrafilter.

\begin{proposition} \label{QmeetsW}
Every $Q$-point is a weak $\W$-ultrafilter.
\end{proposition}

\begin{proof}
Consider partition of $\omega$ into finite sets $\omega = \bigcup_{n \in \omega} I_n$
where $I_0 = \{0, 1\}$ and $I_n=[2^n,2^{n+1})$ for every $n \geq 1$.
Assume $f:\omega \rightarrow \omega$ is an arbitrary finite-to-one function.
If $\U$ is a $Q$-point then there exists $U \in \U$ such that $|U \cap f^{-1}[I_n]| \leq 1$
for every $n \in \omega$. Now enumerate (increasingly) $f[U] = \{u_n: n \in \omega\}$. Either
$U_0 = f^{-1}\{u_{2k}: k \in \omega\}$ or $U_1 = f^{-1}\{u_{2k+1}: k \in \omega\}$ belongs to the ultrafilter $\U$.
Without loss of generality we may assume $U_0 \in \U$. It follows from the definition that $f[U_0]$ does not
contain an arithmetic progression of length 3, thus $f[U_0] \in \W$.
\end{proof}

Weak $\W$-ultrafilter in the previous proposition cannot be replaced with $\W$-ultrafilter because
$Q$-points in general need not be $\W$-ultrafilters. It follows from Proposition 2.4.7 in \cite{FTh} according to which assuming
Martin's axiom for countable posets for every tall ideal $\I$ there exists a $Q$-point which is not an $\I$-ultrafilter.

One can, of course, ask whether the implication in Proposition \ref{QmeetsW} can be reversed. The answer is negative. If Martin's axiom for countable posets holds then
there exists even a $\W$-ultrafilter which is not a $Q$-point. We will construct such an ultrafilter in the rest of this section starting with the following rather technical statement:

\begin{proposition} \label{lemmaWnotQ}
(MA${}_{\hbox{\small ctble}}$) Assume $\F$ is a filter base of cardinality less
than $\ccc$ with the following property
$$(\spadesuit) \qquad (\forall F \in \F) \, (\forall k \in \omega) \, (\exists
n \in \omega) \, |F \cap [2^n, 2^{n+1})| \geq k.$$
Assume $f \in {}^{\omega}\omega$. Then there exists $G \in [\omega]^{\omega}$
such that $f[G] \in \W$ and the filter base generated by $\F$ and $G$ has
property $(\spadesuit)$.
\end{proposition}

\begin{proof}
If there exists $F \in \F$ such that $f[F] \in \W$, let $G=F$.
If there exists $K \in [\omega]^{<\omega}$ such that the filter base
generated by $\F$ and $f^{-1}[K]$ has property $(\spadesuit)$, then put $G = f^{-1}[K]$.
In the following we will assume that no such $K$ exists. From this assumption it follows that
for every $K \in [\omega]^{<\omega}$ and for every $F \in \F$ and for every $k \in \omega$
there exists $n \in \omega$ such that $|\left(F \setminus f^{-1}[K]\right) \cap [2^n, 2^{n+1})|\geq k$.
We construct a suitable $G$ eventually making use of Martin's axiom.
\smallskip

Consider the set $$P = \{K \in [\omega]^{< \omega}: f[K] \hbox{ contains no arithmetic progressions of length }3\}$$
equipped with a partial order $\leq_P$ defined in the following way: $K \leq_P L$ if and only if
$K = L$ or $K \supset L$ and $\min (K \setminus L) > \max L$. For $F \in
\F$ and $k \in \omega$ define $D_{F,k} = \{K \in P: (\exists n \in \omega) \, |K \cap F \cap [2^n, 2^{n+1})| \geq k\}$.
\medskip

{\sl Claim: $D_{F,k}$ is a dense subset of $(P, \leq_P)$ for every $F \in \F$, $k \in \omega$.}
\smallskip

\noindent
{\it Proof of the Claim.\/}
Let $L \in P$ be arbitrary. We want to find $K \in D_{F,k}$ such that $K \leq_P L$.
Let $n_0 = \max \{n \in \omega: L \cap [2^n, 2^{n+1}]) \neq \emptyset\}$.

\medskip
\underline{\it Case I.\/} \ $\sup_{n \in \omega} |f[F \cap [2^n, 2^{n+1})]| = m < \infty$
\smallskip

For $N=\{i \in \omega: i \leq 3\cdot \max f[L]\}$ there exists $n(N) > n_0$ such that $|\left(F \setminus f^{-1}[N]\right) \cap [2^{n(N)}, 2^{n(N)+1})|\geq k\cdot (m+1)$. According to the assumption of Case I. $|f[F \cap [2^{n(N)}, 2^{n(N)+1})]| \leq m$. Now, it follows from Dirichlet's box principle that there exist
$$l \in f[F \cap [2^{n(N)}, 2^{n(N)+1})] \hbox{ and } L' \subseteq \left(F \setminus f^{-1}[N]\right) \cap [2^{n(N)}, 2^{n(N)+1})$$ such that $|L'| \geq k$, $f[L'] =\{l\}$ and $l > 3 \cdot \max f[L]$. Put $K = L \cup L'$. It remains to verify that $K$ is as required:
\begin{enumerate}
\item
Observe that $K \in P$. Any arithmetic progression of length 3 in $f[K]$ contains $l$ (otherwise $f[L]$ contains an arithmetic progression of length 3 in contradiction to the assumption $L \in P$). However, for every $a, b \in f[L]$ we have
$$l - \max\{a, b\} > 2 \cdot \max f[L] > \max f[L] > |a - b|.$$
Thus $a$, $b$ and $l$ cannot form an arithmetic progression.
\item
Obviously, $K \in D_{F,k}$ because $$|K \cap F \cap [2^{n(N)}, 2^{n(N)+1})| \geq |L' \cap F \cap [2^{n(N)}, 2^{n(N)+1})| \geq k.$$
\item
Notice that $K \leq_P L$ because $n(N) > n_0$ and consequently $$\min (K \setminus L) = \min L' > \max L.$$
\end{enumerate}

\smallskip
\underline{\it Case II.\/} \ $\sup_{n \in \omega} |f[F \cap [2^n, 2^{n+1})]| = \infty$
\smallskip

According to the assumption of Case II. there exists $n_1 > n_0$ such that $|f[F \cap [2^{n_1}, 2^{n_1+1})]|\geq 3 \cdot \max f[L] + \left(|L| + k\right)^2$.
Let $$A_0 = \{m \in f[F \cap [2^{n_1}, 2^{n_1+1})]: m > 3 \cdot \max f[L]\}.$$
Obviously, $|A_0| \geq \left(|L|+k\right)^2$.
Now, choose $l_i \in A_0$ for $i = 0, \dots, k-1$ so that $f[L] \cup \{l_i: i < k\}$ contains no arithmetic progressions of length $3$.
This can be done by induction on $i$:

For $i=0$ let $l_0 = \min A_0$. The set $B_1=f[L] \cup \{l_0\}$ does not contain arithmetic progressions of length $3$ because $f[L]$ did not and for arbitrary $a, b \in f[L]$ $$l_0 - \max\{a, b\} > 2 \cdot \max f[L] > \max f[L] > |a - b|.$$

If $0<i<k-1$ and $l_j \in A_0$ for $j < i$ are already known such that $B_i = f[L] \cup \{l_j: j < i\}$ contains no arithmetic progressions of length $3$, define
$$A_{i}=\{m \in A_0: (\exists a,b \in B_i) \ a, b, m \hbox{ form an arithmetic progression}\}.$$
Since $|B_i| = |L|+i$ the set $A_i$ has at most $\frac{1}{2}(|L|+i)(|L|+i-1)$ elements.
So $A_0 \setminus A_i \neq \emptyset$ and we may define $l_i = \min (A_0 \setminus A_i)$. It follows from the construction that $B_{i+1}=B_i \cup \{l_i\}$ contains no arithmetic progressions of length $3$.

Finally, let $L' = F \cap [2^{n_1}, 2^{n_1+1}) \cap f^{-1}[\{l_i: i <k\}]$ and put $K = L \cup L'$.
It remains to verify that $K$ is as required:
\begin{enumerate}
\item
Obviously, $K \in P$ because $f[K] = f[L] \cup \{l_i: i < k\}$ and the latter set contains no arithmetic progressions of length 3.
\item
Observe that $K \in D_{F,k}$ because $$|K \cap F \cap [2^{n_1}, 2^{n_1+1})| \geq |L'| \geq k.$$
\item
Notice that $K \leq_P L$ because $n_1 > n_0$ and consequently
\end{enumerate}
{\hfill $\min (K \setminus L) = \min L' > \max L.$ \hfill $\Box$ {\small \sl Claim.\/}}

\bigskip

Since the family $\D = \{D_{F,k}: F \in \F, k \in \omega\}$ consists of dense
subsets of the countable poset $P$ and $|\D| < \ccc$, it follows from Martin's
axiom for countable posets that there is a $\D$-generic filter $\G$.

Let $G = \bigcup\{K: K \in \G\}$. It remains to verify that $G$ is as required, i.e.
\begin{enumerate}[a)]
\item $f[G]$ contains no arithmetic progressions of length 3, thus $f[G] \in \W$
\item $(\forall F \in \F) \, (\forall k \in \omega) \, (\exists
n \in \omega) \, |F \cap G \cap [2^n, 2^{n+1})| \geq k.$
\end{enumerate}

\noindent
For a): Consider $a, b, c \in f[G]$ arbitrary. There exist $K_a$, $K_b$, $K_c \in \G$ such that $a \in f[K_a]$, $b \in f[K_b]$ and $c \in f[K_c]$.
Since $\G$ is a filter, there exists $K_0 \in \G$ which is $\leq_P$-below all three sets  $K_a$, $K_b$, $K_c$. Thus $a, b, c \in f[K_0]$.
Because $K_0$ is an element of $P$, the set $f[K_0]$ contains no arithmetic progressions of length 3. In particular $a$, $b$ and $c$ do not form an arithmetic progression.

\medskip

\noindent
For b): Take $k \geq 1$ arbitrary. For every $K \in \G \cap D_{F,k}$ we have $G
\supseteq K$ and $ |F \cap G \cap [2^n, 2^{n+1})| \geq |F \cap K \cap [2^n, 2^{n+1})| \geq k$ for some $n$.
\end{proof}

\begin{lemma} \label{extNotQ}
Every filter base $\F$ which has property $(\spadesuit)$ introduced in Proposition \ref{lemmaWnotQ}
can be extended into an ultrafilter which is not a $Q$-point.
\end{lemma}

\begin{proof}
The family $\{[2^n, 2^{n+1}): n \in \omega\}$ is a partition of $\omega$ into finite sets witnessing the fact that an ultrafilter with property $(\spadesuit)$ is not a $Q$-point. Therefore we will show that every filter base $\F$ with property $(\spadesuit)$ can be extended into an ultrafilter with property $(\spadesuit)$.
This can be accomplished by transfinite induction on $\alpha < \ccc$ where in each non-limit step one subset of $\omega$ is considered and the filter base is extended by either the set itself or its complement.

To this end, consider a filter base $\F$ with property $(\spadesuit)$ and $A \subseteq \omega$: Either $(\forall F \in \F) \, (\forall k \in \omega) \, (\exists n \in \omega) \, |F \cap A \cap [2^n, 2^{n+1})| \geq k$ and then the filter base $\F'$ generated by $\F$ and $A$ has property $(\spadesuit)$.\\
Or $(\exists F_0 \in \F) \, (\exists k_0 \in \omega) \, (\forall n \in \omega) \, |F_0 \cap A \cap [2^n, 2^{n+1})| < k_0$ and then the filter base $\F'$ generated by $\F$ and $\omega \setminus A$ has property $(\spadesuit)$. Indeed, since $\F$ has property $(\spadesuit)$ for every $F \in \F$ and for every $k \in \omega$ there exists $n \in \omega$ such that $|F \cap F_0 \cap [2^n, 2^{n+1})| \geq k+k_0$. However, $|F \cap F_0 \cap A \cap [2^n, 2^{n+1})| < k_0$ and so $|F \cap (\omega \setminus A) \cap [2^n, 2^{n+1})| \geq |F \cap F_0 \cap (\omega \setminus A) \cap [2^n, 2^{n+1})| \geq k$.
\end{proof}

\begin{theorem} \label{thmWnotQ}
(MA${}_{\hbox{\small ctble}}$) There is a $\W$-ultrafilter which is not a $Q$-point.
\end{theorem}

\begin{proof}
Enumerate ${}^{\omega}\omega$ as $\{f_{\alpha}:\alpha < \ccc\}$.
By transfinite induction on $\alpha < \ccc$ we will construct filter bases
$\F_{\alpha}$ so that the following conditions are satisfied:
\smallskip

(i) $\F_0$ is the Fr\'echet filter

(ii) $\F_{\alpha} \subseteq \F_{\beta}$ whenever $\alpha \leq \beta$

(iii) $\F_{\gamma} = \bigcup_{\alpha < \gamma} \F_{\alpha}$ for $\gamma$ limit

(iv) $(\forall \alpha)$ $|\F_{\alpha}| \leq (|\alpha| + 1) \cdot \omega$

(v) $(\forall \alpha)$ $\F_{\alpha}$ has property $(\spadesuit)$

(vi) $(\forall \alpha)$ $(\exists F \in \F_{\alpha+1})$ $f_{\alpha}[F] \in \W$
\smallskip

Condition (i) starts the induction and (ii), (iv) allow it to keep going.
Limit stages are taken care of by condition (iii) so it remains to show that successor stages can be handled.
\smallskip

Successor stage: Suppose we already know $\F_{\alpha}$. If there is $A \in
\F_{\alpha}$ such that $f_{\alpha}[A] \in \W$ then simply put
$\F_{\alpha + 1} = \F_{\alpha}$.  If $f_{\alpha}[F] \not\in \W$
for every $F \in \F_{\alpha}$ then apply Proposition \ref{lemmaWnotQ} to the filter base
$\F_{\alpha}$ and the function $f_{\alpha}$.
Let $\F_{\alpha+1}$ be the filter base generated by $\F_{\alpha}$ and $G$.

Finally, let $\F = \bigcup_{\alpha < \ccc} \F_{\alpha}$.
Since the filter base $\F$ has property $(\spadesuit)$, it can be extended into an
ultrafilter which is not a $Q$-point according to Lemma \ref{extNotQ}.
It is a $\W$-ultrafilter which is not a $Q$-point because each ultrafilter
which extends $\F$ is a $\W$-ultrafilter due to condition (vi).
\end{proof}

\section{Rapid ultrafilters and the van der Waerden ideal}

In this section we prove that rapid ultrafilters, unlike $Q$-points, may have an empty intersection
with the van der Waerden ideal. We will construct such an ultrafilter assuming Martin's axiom for countable posets
in Theorem \ref{MainThm}, which is actually slightly stronger.
\smallskip

Let us start with the definition of summable ideals. They play an important role in an alternative
characterization of rapid ultrafilters, which we use later in the proof.\\
For a function $g: \omega \rightarrow (0,\infty)$ such that $\sum\limits_{n \in \omega} g(n) = +\infty$ the family
$$\I_g = \{A \subseteq \omega: \sum_{a\in A} g(a) < +\infty\}$$
is a {\it summable ideal\/} on $\omega$ {\it determined by function $g$\/}. A summable ideal $\I_g$ is tall if and only if $\lim\limits_{n \rightarrow \infty} g(n) = 0$.

A characterization of rapid ultrafilters involving summable ultrafilters can be found e.g. in \cite{M-A} as Theorem 2.8.10.
We restate it here in terms of weak $\I$-ultrafilters.

\begin{proposition} [\cite{M-A}] \label{rapDescr}
The following statements are equivalent for an ultrafilter $\U \in {\omega}^{\ast}$:
\begin{enumerate}
\item $\U$ is a rapid ultrafilter
\item $\U$ is a weak $\I_g$-ultrafilter for every tall summable ideal $\I_g$
\item $\U$ is a weak $\I$-ultrafilter for every tall analytic $P$-ideal $\I$
\end{enumerate}
\end{proposition}

Notice that analytic ideals in the third clause of the previous proposition are $P$-ideals.
Since the van der Waerden ideal is not a $P$-ideal, the theorem does not contradict
the existence of rapid ultrafilters, whose intersection with the van der Waerden ideal is empty.

Rapid ultrafilters which are disjoint with the van der Waerden ideal, i.e. rapid ultrafilters
which contain only AP-sets, do actually exist if we assume Martin's axiom for countable posets.
Under the same assumption even more is true and there exist hereditarily rapid ultrafilters consisting of AP-sets.
Hereditarily rapid ultrafilters, as the name suggests, form a subclass of rapid ultrafilters, which we define as follows:\\
An ultrafilter $\U$ is a {\it hereditarily rapid
ultrafilter\/} if it is a rapid ultrafilter such that the ultrafilter
$f(\U)$ generated by the sets $\{f[U]: U \in \U\}$ is again a rapid
ultrafilter for every $f \in {}^{\omega}\omega$.\\
The characterization of rapid ultrafilters in Proposition \ref{rapDescr} can be easily reformulated to hereditarily rapid ultrafilters:

\begin{proposition} \label{herRapDescr}
The following statements are equivalent for $\U \in {\omega}^{\ast}$:
\begin{enumerate}
\item $\U$ is a hereditarily rapid ultrafilter
\item $\U$ is an $\I_g$-ultrafilter for every tall summable ideal $\I_g$
\item $\U$ is an $\I$-ultrafilter for every tall analytic $P$-ideal $\I$
\end{enumerate}
\end{proposition}

We will concentrate on the first two clauses of the previous proposition and we will construct our ultimate goal -- a hereditarily
rapid ultrafilter which does not intersect the van der Waerden ideal -- as an ultrafilter
which is an $\I_g$-ultrafilter for every tall summable ideal $\I_g$. The following lemma is crucial
for the successor stages of the construction:

\begin{lemma} \label{lemmaIgnotW}
(MA${}_{\hbox{\small ctble}}$) Let $\I_g$ be a tall summable ideal. Assume $\F$
is a filter base, $|\F|< \ccc$, $\F \cap \W = \emptyset$ and $f \in
{}^{\omega}\omega$.  Then there exists $G \in [\omega]^{\omega}$ such that
$f[G] \in \I_g$ and $F \cap G$ is an AP-set for every $F \in \F$.
\end{lemma}

\begin{proof}
If there exists $F \in \F$ such that $f[F] \in \I_g$, put $G = F$.
If there exists $K \in [\omega]^{< \omega}$ such that $F \cap f^{-1}[K]$ is an AP-set for every $F \in
\F$, then put $G = f^{-1}[K]$.
In the following we will assume that no such set exists, i.e.
$$(\clubsuit) \hbox{ for every } K \in [\omega]^{< \omega} \hbox{ there is } F_K \in \F \hbox{ such that } F_K \cap f^{-1}[K] \in \W.$$
This also means that $F \setminus f^{-1}[K]$ is an AP-set for every $F \in \F$.
\smallskip

Consider the set $$P = \{K \in [\omega]^{<\omega} : \sum\limits_{a \in f[K]} g(a) \leq \left(2 - \frac{1}{2^{|K|}}\right) \cdot \max\limits_{a \in f[K]} g(a)\}$$
with a partial order $\leq_P$ defined by: $K \leq_P L$ if and only if $K = L$ or $K \supset L$ and $\min (K \setminus L) > \max L$.
Now, for every $F \in \F$ and every $k \geq 1$ define $D_{F,k} = \{K \in P : K \cap F$ contains an arithmetic progression of length $k\}$.
\medskip

{\sl Claim. $D_{F,k}$ is dense in $P$ for every $F \in \F$ and $k \geq 1$.}
\smallskip

\noindent
{\it Proof of the Claim.\/}
Take $L \in P$ arbitrary. Since $\lim_{n \rightarrow \infty} g(n)= 0$ there exists
$n_L \in \omega$ such that for every $n > n_L$
$$g(n) < \frac{1}{2^{|L|+1} \cdot k} \cdot \max_{a \in f[L]} g(a).$$

According to the assumption $(\clubsuit)$ there exists $F_{n_L} \in \F$ such that
$F_{n_L} \cap f^{-1}[0, n_L] \in \W$.  It follows that $A_{n_L} = (F \cap
F_{n_L}) \setminus f^{-1}[0, n_L]$ is an infinite AP-set, thus one can choose
an arithmetic progression $L' \subset A_{n_L}$ such that $|L'|=k$ and $\min L'
> \max L$. Let $K = L \cup L'$. Observe that due to the choice of $L'$ one has
$\max\limits_{a \in f[K]} g(a) = \max\limits_{a \in f[L]} g(a)$.
To see that $K \in P$ notice that

$$\sum\limits_{a \in f[K]} g(a) \leq \sum\limits_{a \in f[L]} g(a) +
\sum\limits_{a \in f[L']} g(a) \leq$$
$$\leq \left(2 - \frac{1}{2^{|L|}}\right) \cdot \max\limits_{a \in f[L]} g(a) +
|L'| \cdot \max\limits_{a \in f[L']} g(a) \leq$$
$$\leq \left(2 - \frac{1}{2^{|L|}}\right) \cdot \max\limits_{a \in f[L]} g(a) +
k \cdot \frac{1}{2^{|L|+1} \cdot k} \cdot \max\limits_{a \in f[L]} g(a) \leq$$
$$\leq \left(2 - \frac{1}{2^{|L|+1}}\right)\max\limits_{a \in f[L]} g(a)
\leq \left(2 - \frac{1}{2^{|K|}}\right)\max\limits_{a \in f[K]} g(a)$$

It is obvious that $K \leq_P L$. Also $K \in D_{F,k}$ because $K \cap F \supseteq L' \cap F$
contains an arithmetic progression of length $k$. \hfill $\Box$ {\small \sl Claim.}
\medskip

Since the family $\D = \{D_{F,k}: F \in \F, k \geq 1\}$ consists of dense
subsets of the countable poset $P$ and $|\D| < \ccc$, it follows from Martin's
axiom for countable posets that there is a $\D$-generic filter $\G$. 

Let $G = \bigcup\{K: K \in \G\}$. It remains to verify that $G$ is as required:
\begin{enumerate}[a)]
\item $f[G] \in \I_g$, i.e.~$\sum_{a \in f[G]} g(a) < +\infty$
\item $(\forall F \in \F)$ $G \cap F$ is an AP-set
\end{enumerate}

\noindent
For a): Enumerate $f[G] = \{u_n: n \in \omega\}$. For every $n$ there exists $K_n \in
\G$ such that $u_n \in f[K_n]$. We may assume $K_{n+1} \leq_P K_n$ (and thus
$f[K_{n+1}] \supseteq f[K_n]$) for every $n \in \omega$ because $\G$ is a
filter.  Since $f[G] = \bigcup_{n \in \omega} f[K_n]$, we get
$$\sum\limits_{a \in f[G]} g(a) = \lim\limits_{n \rightarrow \infty} \sum\limits_{a \in f[K_n]}
g(a) \leq \lim\limits_{n \rightarrow \infty}(2 - \frac{1}{2^{|K_n|}})
\max\limits_{a \in f[K_n]} g(a) \leq 2 \cdot \max\limits_{a \in \omega} g(a).$$

\medskip

\noindent
For b): Take $k \geq 1$ arbitrary. For every $K \in \G \cap D_{F,k}$ we have $G
\supseteq K$ and $K \cap F$ contains an arithmetic progression of length $k$.
Hence $G \cap F$ contains arithmetic progressions of arbitrary length, i.e.~$G \cap F$ is an AP-set.
\end{proof}
						
\begin{theorem}  \label{MainThm}
(MA${}_{\hbox{\small ctble}}$) There is a hereditarily rapid ultrafilter $\U$ such that
$\U \cap \W = \emptyset$.
\end{theorem}
						
\begin{proof}
Enumerate as  $\{\langle f_{\alpha}, g_{\alpha} \rangle : \alpha <
\ccc\}$ all pairs $\langle f_{\alpha}, g_{\alpha} \rangle$ where
$f_{\alpha} \in {}^{\omega}\omega$ and $\I_{g_{\alpha}}$ is a
tall summable ideal. By transfinite induction on $\alpha < \ccc$ we will
construct filter bases $\F_{\alpha}$ so that the following conditions are
satisfied:
\smallskip

(i) $\F_0$ is the Fr\'echet filter

(ii) $\F_{\alpha} \subseteq \F_{\beta}$ whenever $\alpha \leq \beta$

(iii) $\F_{\gamma} = \bigcup_{\alpha < \gamma} \F_{\alpha}$ for $\gamma$ limit

(iv) $(\forall \alpha)$ $|\F_{\alpha}| \leq (|\alpha| + 1) \cdot \omega$

(v) $(\forall \alpha)$ $(\forall F \in \F_{\alpha})$ $F$ is an AP-set

(vi) $(\forall \alpha)$ $(\exists F \in \F_{\alpha+1})$ $f_{\alpha}[F] \in \I_{g_{\alpha}}$
\smallskip

Condition (i) starts the induction and (ii), (iv) allow it to keep going.
Limit stages are taken care of by condition (iii) so it remains to show that successor stages can be handled.
\smallskip

Successor stage: Suppose we already know $\F_{\alpha}$. If there is $A \in
\F_{\alpha}$ such that $f_{\alpha}[A] \in \I_{g_{\alpha}}$ then simply put
$\F_{\alpha + 1} = \F_{\alpha}$.  If $f_{\alpha}[F] \not\in \I_{g_{\alpha}}$
for every $F \in \F_{\alpha}$ then apply Lemma \ref{lemmaIgnotW} to the ideal
$\I_{g_{\alpha}}$, the filter base $\F_{\alpha}$ and the function $f_{\alpha}$.
Let $\F_{\alpha+1}$ be the filter base generated by $\F_{\alpha}$ and $G$.

Finally, let $\F = \bigcup_{\alpha < \ccc} \F_{\alpha}$. Since $F$ is an AP-set
for every $F \in \F$ the filter base $\F$ can be extended to an ultrafilter $\U$
which does not intersect the van der Waerden ideal. Every ultrafilter
which extends $\F$, however, is an $\I_{g_{\alpha}}$-ultrafilter for every tall
summable ideal $\I_{g_{\alpha}}$ because of condition (vi). Thus $\U$ is
a hereditarily rapid ultrafilter satisfying $\U \cap \W = \emptyset$.
\end{proof}
												
\begin{corollary}
(MA${}_{\hbox{\small ctble}}$) There exists a rapid ultrafilter $\U$ such that
$\U$ contains only AP-sets. \hfill $\Box$
\end{corollary}

\end{document}